\documentclass[11pt, reqno]{amsart}
\usepackage{amsmath,amssymb,amsfonts,amscd,hyperref,color}
\usepackage[latin1]{inputenc}

\usepackage{hyperref}
\hypersetup{
     colorlinks   = true,
     citecolor    = blue
}

\usepackage{mathtools}
\DeclarePairedDelimiter\ceil{\lceil}{\rceil}

\usepackage[abs]{overpic}

\usepackage{url}

\usepackage{verbatim}

\newcommand{\Hmm}[1]{\leavevmode{\marginpar{\tiny%
$\hbox to 0mm{\hspace*{-0.5mm}$\leftarrow$\hss}%
\vcenter{\vrule depth 0.1mm height 0.1mm width \the\marginparwidth}%
\hbox to
0mm{\hss$\rightarrow$\hspace*{-0.5mm}}$\\\relax\raggedright #1}}}

\newcommand\blue[1]{\textcolor{blue}{#1}}

\newtheorem{thm}{Theorem}[section]
\newtheorem{cor}[thm]{Corollary}

\newtheorem{lem}[thm]{Lemma}

\theoremstyle{definition}
\newtheorem{defi}[thm]{Definition}

\newtheorem{rem}[thm]{Remark}

\numberwithin{equation}{section}

\begin{document}

\title{Discrete Weighted Hardy Inequality in 1-D}
\markright{Discrete Weighted Hardy Inequality}

\author[Shubham Gupta]{By Shubham Gupta}
\address{(Shubham Gupta) Department of Mathematics, Imperial College London, 180 Queen's Gate, London, SW7 2AZ, United Kingdom.}
\email{s.gupta19@imperial.ac.uk}

\maketitle
\hypersetup{linkcolor=blue}

\vspace{-30pt}
\begin{abstract}
In this paper we consider weighted versions of one dimensional discrete Hardy's inequality on the half-line with \emph{power weights} of the form $n^\alpha$; namely, we consider:
\begin{equation}\label{0.1}
    \sum_{n=1}^\infty |u(n)-u(n-1)|^2 n^\alpha \geq c(\alpha) \sum_{n=1}^\infty \frac{|u(n)|^2}{n^2}n^\alpha.
\end{equation}

We prove the above inequality when $\alpha \in [0,1) \cup [5,\infty)$ with the sharp constant $c(\alpha)$. Furthermore, when $\alpha \in [1/3,1) \cup \{0\}$, we prove an improved version of \eqref{0.1} by adding infinitely many positive lower order terms in the RHS. More precisely, we prove 
\begin{equation}
    \sum_{n=1}^\infty |u(n)-u(n-1)|^2 n^\alpha \geq c(\alpha) \sum_{n=1}^\infty \frac{|u(n)|^2}{n^2} n^\alpha + \sum_{k=3}^\infty b_k(\alpha) \sum_{n=2}^\infty \frac{|u(n)|^2}{n^k}n^\alpha
\end{equation}
for non-negative constants $b_k(\alpha)$.

\end{abstract}

\section{Introduction}\label{sec1}

\let\thefootnote\relax\footnotetext{Keywords: Hardy's inequality, Super-solution method, Sharp constant, Power weights.\vspace{3pt}

2020 Mathematics Subject Classification: 39B62, 26D15.}

In 1921 Landau wrote a letter to G.H. Hardy including proof of the following inequality with the sharp constant\cite{landau1921letter}: 
\begin{equation}\label{1.1}
    \sum_{n=1}^\infty a_n^p \geq 
    \Big( \frac{p-1}{p}\Big)^p \sum_{n=1}^\infty \Big(\frac{a_1+a_2+...+a_n}{n}\Big)^p  
\end{equation}
for $p>1$ where $\{a_n\}_{n=1}^\infty$ is an arbitrary non-negative sequence of real numbers. \\

This inequality is referred to as Hardy's inequality since then (see \cite{kufner2006prehistory} for a beautiful description of the prehistory of Hardy's Inequality). The author would also like to mention a recent and short proof of \eqref{1.1} by Lefevre \cite{lefevre2020}.\\

Let $C_c(\mathbb{N}_{0})$ be the space of finitely supported functions on $\mathbb{N}_0 = \{0,1,2,3,...\}$. It is not very hard to see that for $p>1$, \eqref{1.1} is equivalent to
\begin{equation}\label{1.2}
    \sum_{n=1}^\infty|u(n)-u(n-1)|^p \geq
    \Big( \frac{p-1}{p}\Big)^p \sum_{n=1}^\infty \frac{|u(n)|^p}{|n|^p}
\end{equation}
for all $u \in C_c(\mathbb{N}_0)$ with the ``Dirichlet Boundary Condition" $u(0)=0$. Recently \eqref{1.2} was improved for the case $p=2$ \cite{pinchover2018improved}, and later for general $p>1$ in \cite{keller2019improved}. More precisely, authors in \cite{keller2019improved} prove the following result:
\begin{equation}\label{1.3}
    \sum_{n=1}^\infty |u(n)-u(n-1)|^2 \geq \frac{1}{4}\sum_{n=1}^\infty \frac{|u(n)|^2}{n^2} + \sum_{k=2}^\infty{4k \choose 2k} \frac{1}{(4k-1)2^{4k-1}} \sum_{n=2}^\infty \frac{|u(n)|^2}{n^{2k}}.
\end{equation}

Although there is extensive literature on the continuous analogues of Hardy's inequality \eqref{1.2}(see classical books \cite{belbook, kufnerbook, Mazyabook} and references therein), very little is known about these inequalities in the discrete setting. One of the major hurdles is that calculus breaks down in the discrete setting, making it difficult to extend proofs of Hardy's inequality in the continuum to the discrete setting. It is worthwhile to mention the works \cite{laptev2016}, \cite{pinchover2018}, \cite{berchio2021poincare} which has been successful in overcoming the absence of calculus. In \cite{laptev2016} Kapitanski and Laptev studied discrete Hardy's inequality of the form \eqref{1.2} on higher dimensional grids  $\mathbb{Z}^d$ by converting it to a problem on the torus using Fourier transform methods. In \cite{pinchover2018} Keller et al. proved Hardy-type inequalities on general graphs with \emph{optimal} weights by developing a discrete version of the \emph{super-solution} method. Recently the method used in \cite{pinchover2018} was exploited to prove some new discrete Hardy's inequalities on regular trees in \cite{berchio2021poincare}. Before getting into the main setting of the paper, we would like to quote papers \cite{gao2012}, \cite{liu2012}, \cite{braverman1994}, \cite{miclo1999}, \cite{bui2020application}, \cite{keller2021hardy}, \cite{keller2021optimal}, \cite{kostenko2021heat}  where various variants of \eqref{1.2} are considered, improved and  applied. \\ 

The goal of this paper is to prove weighted versions of inequality \eqref{1.2} and \eqref{1.3} for the case $p=2$ with \emph{power weights} $n^\alpha$. One of the main results of this paper is the following two-parameter family of weighted Hardy's inequalities: If $\alpha, \beta \in \mathbb{R}$ then
\begin{equation}\label{1.4}
    \sum_{n=1}^\infty |u(n)-u(n-1)|^2 n^\alpha \geq \sum_{n=1}^\infty w_{\alpha,\beta}(n) |u(n)|^2
\end{equation}
where 
\begin{equation}\label{1.5}
    w_{\alpha, \beta}(n) := n^\alpha \Bigg[ 1 + \Big(1+\frac{1}{n}\Big)^\alpha - \Big(1-\frac{1}{n}\Big)^\beta - \Big(1+\frac{1}{n}\Big)^{\alpha+\beta}\Bigg]
\end{equation}
for $n \geq 2$ and $w_{\alpha,\beta}(1) := 1+2^\alpha - 2^{\alpha+\beta}$. \\

As will be shown, \eqref{1.4} contains the following power weights Hardy Inequalities as special cases:
\begin{equation}\label{1.6}
    \sum_{n=1}^\infty |u(n)-u(n-1)|^2 n^\alpha \geq \frac{(\alpha-1)^2}{4}\sum_{n=1}^\infty \frac{|u(n)|^2}{n^2} n^\alpha
\end{equation}
whenever $\alpha \in[0,1)$ or $ \alpha \in [5,\infty)$.\\
and we have an improved version of \eqref{1.6} for $\alpha \in [1/3,1) \cup \{0\}$
\begin{equation}\label{1.7}
    \sum_{n=1}^\infty |u(n)-u(n-1)|^2 n^\alpha \geq \frac{(\alpha-1)^2}{4} \sum_{n=1}^\infty \frac{|u(n)|^2}{n^2} n^\alpha + \sum_{k=3}^\infty b_k(\alpha) \sum_{n=2}^\infty \frac{|u(n)|^2}{n^k}n^\alpha
\end{equation}
where the non-negative constants $b_k(\alpha)$ are given by
\begin{equation}\label{1.8}
    b_k(\alpha) := {\alpha \choose k} - (-1)^k {(1-\alpha)/2 \choose k} -{(1+\alpha)/2 \choose k}.
\end{equation}

\begin{rem}\label{rem1.1}
Inequality \eqref{1.6} is derived from \eqref{1.4} by estimating $w_{\alpha,\beta}$ by $\frac{(\alpha-1)^2}{4} n^{\alpha-2}$ from below by choosing $\beta = (1-\alpha)/2$. We would like to point out that this lower estimate on $w_{\alpha, \beta}$ fails to hold true when $\alpha <0$ or $\alpha \in(1,4)$(this will be proved in section \ref{sec5} of the paper). Due to this reason we fail to prove \eqref{1.6} for all non-negative $\alpha$. With the aim of proving inequalities of type \eqref{1.6} for all $\alpha \geq 0$, one could ask the following question: Is it possible to find $\beta$ and non-negative constant $c(\alpha)$ such that $w_{\alpha, \beta}(n) \geq c(\alpha) n^{\alpha-2}$? We couldn't manage to answer this question in this paper.      
\end{rem}

\begin{rem}\label{rem1.2}
We would like to mention that \eqref{1.7} is true for all $\alpha \in  [0,1) \cup [5, \infty)$ but we conjecture that the constant $b_k(\alpha)$ is not non-negative for all $k \geq 3$ when $\alpha$ lies outside $[1/3,1) \cup\{0\}$, that is, when $\alpha \in (0,1/3) \cup [5, \infty)$(it will be partially proved in section \ref{sec5}).
\end{rem}

Our approach is based on the \emph{supersolution} method. This is a well known method for proving Hardy-type Inequalities in the continuous setting.(see \cite{cazacu2020} for the survey of methods known for proving Hardy-type inequalities in the continuum). Let us sketch briefly the idea behind the supersolution method. The standard Hardy-Inequality in the continuous setting states 
\begin{equation}\label{1.9}
    \int_{\mathbb{R}^d} |\nabla u|^2 dx \geq \frac{(d-2)^2}{4} \int_{\mathbb{R}^d} \frac{|u(x)|^2}{|x|^2} dx. 
\end{equation}
for all $u \in C_c^\infty(\mathbb{R}^d)$ and $d \geq 3$. The super-solution method to prove \eqref{1.9} roughly goes as follows. Let $u = \varphi \psi$. Then 
\begin{align*}
    |\nabla u|^2 = \psi^2 |\nabla \varphi|^2 + \varphi^2 |\nabla \psi|^2 + 2 \nabla \varphi \cdot \nabla \psi \varphi \psi.
\end{align*}
Applying integration by parts we obtain
\begin{align*}
    \int |\nabla u|^2 & = \int \psi^2 |\nabla\varphi|^2 + \int \phi^2 |\nabla \psi|^2 + 1/2\int \nabla(\varphi^2) \cdot \nabla (\psi^2)  \\
    &= \int \varphi^2 |\nabla \psi|^2 - \int \varphi \psi^2 \Delta \varphi  \geq \int \frac{-\Delta \varphi}{\varphi}|u|^2. 
\end{align*}
If $\varphi$ satisfies $\frac{-\Delta \varphi}{\varphi} \geq w$ then we have
\begin{equation}\label{1.10}
    \int |\nabla u|^2 dx \geq \int  w(x) |u|^2 dx.
\end{equation}
Therefore proving \eqref{1.7} boils down to a much simpler task of finding a solution of $ - \Delta \varphi - w \varphi \geq 0$ with $w = \frac{c}{|x|^2}$. This simple idea of connecting Hardy-type inequalities with solution of differential equations has been exploited a lot in the literature to prove various weighted version and improvements of first-order inequalities of the form \eqref{1.9}(\cite{ghoussoub2011}, \cite{Ghoussoub2013}). In this paper we prove a discrete version of the supersolution method which will then be used to prove one of the main result \eqref{1.4}. \\\\
The paper is divided into various sections. In section \ref{sec2} we will properly state the main results of the paper. In section \ref{sec3} we derive the discrete analogue of supersolution method and using that we will prove \eqref{1.4}. In section \ref{sec4} we derive the inequalities \eqref{1.6} and \eqref{1.7} from the \eqref{1.4}. Finally in section \ref{sec5} we will comment a bit about the limitation of the method: proving the results mentioned in the remarks \ref{rem1.1} and \ref{rem1.2}. \\

\textbf{Acknowledgements.} 
I am grateful to Professor Ari Laptev for suggesting the problem and for various valuable discussions. I would also like to thank him for comments on the early drafts of this paper. Finally, I thank the reviewers for their thorough reading and many helpful suggestions. The author is funded by President's Ph.D. Scholarship, Imperial College London.

\section{Main Results}\label{sec2}
The first main result is the following two-parameter family of discrete weighted Hardy inequalities.
\begin{thm}\label{thm2.1}
If $\alpha, \beta \in \mathbb{R}$, then
\begin{equation}\label{2.1}
    \sum_{n=1}^\infty |u(n)-u(n-1)|^2 n^\alpha \geq \sum_{n=1}^\infty w_{\alpha,\beta}(n) |u(n)|^2
\end{equation}
for $u \in C_c(\mathbb{N}_0)$ and $u(0)=0$, \\ 
where 
\begin{equation}\label{2.2}
    w_{\alpha, \beta}(n) := n^\alpha \Bigg[ 1 + \Big(1+\frac{1}{n}\Big)^\alpha - \Big(1-\frac{1}{n}\Big)^\beta - \Big(1+\frac{1}{n}\Big)^{\alpha+\beta}\Bigg]
\end{equation}
for $n \geq 2$ and $w_{\alpha, \beta}(1) := 1+2^\alpha - 2^{\alpha+\beta}$.
\end{thm}
\begin{rem}
We would like to mention that inequality \eqref{2.1} is a generalization of improved Hardy's inequality in \cite{pinchover2018improved}. We recover the inequality in \cite{pinchover2018improved}, by taking $\alpha=0$ and $\beta =1/2$ in inequality \eqref{2.1}. 
\end{rem}

As a special case of Theorem \ref{thm2.1}, we obtain the following power weight discrete Hardy's Inequality:
\begin{cor}\label{cor2.3}
Let $\alpha \in [0,1) \cup [5,\infty)$. Then for all $u \in C_c(\mathbb{N}_0)$ with $u(0)=0$ we have \begin{equation}\label{2.3}
    \sum_{n=1}^\infty |u(n)-u(n-1)|^2 n^\alpha \geq \frac{(\alpha-1)^2}{4} \sum_{n=1}^\infty \frac{|u(n)|^2}{n^2} n^\alpha
\end{equation}
Moreover the constant in \eqref{2.3} is sharp; that is, if we replace $(\alpha-1)^2/4$ with a strictly bigger constant then inequality \eqref{2.3} will not be true. 
\end{cor}
\begin{rem}
Note that inequality \eqref{2.3} with $\alpha =0$ yields classical discrete Hardy's inequality \eqref{1.2} for $p=2$.
\end{rem}

Inequality \eqref{2.1} also yields the following improvement of \eqref{2.3} when $\alpha \in [1/3,1)\cup \{0\}$.
\begin{cor}\label{cor2.5}
If $\alpha \in [1/3,1) \cup\{0\}$ then 
\begin{equation}\label{2.4}
    \sum_{n=1}^\infty |u(n)-u(n-1)|^2 n^\alpha \geq \frac{(\alpha-1)^2}{4} \sum_{n=1}^\infty \frac{|u(n)|^2}{n^2} n^\alpha + \sum_{k=3}^\infty b_k(\alpha) \sum_{n=2}^\infty \frac{|u(n)|^2}{n^k}n^\alpha 
\end{equation}
for all $u \in C_c(\mathbb{N}_0)$ with $u(0)=0$,\\ 
where the non-negative coefficients $b_k(\alpha)$ are given by
\begin{equation}\label{2.5}
    b_k(\alpha) := {\alpha \choose k} - (-1)^k {(1-\alpha)/2 \choose k} -{(1+\alpha)/2 \choose k}
\end{equation}
where ${\gamma \choose r}$ is the binomial coefficient for real parameters $\gamma$ and $r$.
\end{cor}
\begin{rem}
Inequality \eqref{2.4} for $\alpha =0$ follows from the improved Hardy inequality proved in \cite{pinchover2018improved}. In fact inequality proved in \cite{pinchover2018improved} is strictly stronger than \eqref{2.4} for $\alpha =0$.     
\end{rem}

\section{Discrete Super-Solution Method}\label{sec3}
\begin{defi}\label{ef}
Let $\varphi$ be a real-valued function on $\mathbb{N}_0$. Then the \emph{combinatorial laplacian} $\Delta$ is defined as 
\begin{align*}
    \Delta \varphi(n):= 
    \begin{cases}
        \varphi(n) - \varphi(n-1) + \varphi(n)-\varphi(n+1) \hspace{19pt} \text{for } \hspace{5pt}  n \geq 1\\
        \varphi(n)-\varphi(n+1) \hspace{109pt} \text{for } \hspace{5pt}  n=0
    \end{cases}
\end{align*}
\end{defi}

\begin{lem} \label{lem3.1} Let $v$ and $w$ be non-negative functions on $\mathbb{N}$. Assume $\exists$ function $\varphi : \mathbb{N}_0 \rightarrow [0, \infty)$ which is positive on $\mathbb{N}$ such that
\begin{equation}\label{3.1}
    \Big(\Delta \varphi(n) v(n) - (\varphi(n+1)-\varphi(n))(v(n+1)-v(n))\Big) \geq w(n)\varphi(n) 
\end{equation}
for all $n \in \mathbb{N}$. Then following inequality holds true
\begin{equation}\label{3.2}
    \sum_{n=1}^\infty |u(n)-u(n-1)|^2 v(n) \geq \sum_{n=1}^\infty w(n)|u(n)|^2
\end{equation}
for $ u \in C_c(\mathbb{N}_0)$ and $u(0)= 0$.
\end{lem}

\begin{proof}
It can be easily seen that for $a \in \mathbb{R}$ and $t \geq 0$ we have
\begin{equation}\label{3.3}
    (a-t)^2 \geq (1-t)(a^2 -t).
\end{equation}
Let $\psi(n) := \frac{u(n)}{\varphi(n)} $ on $\mathbb{N}$ and $\psi(0):=0$. Assuming $\psi(m) \neq 0$ and applying \eqref{3.3} for $a = \psi(n)/\psi(m)$ and $t = \varphi(m)/\varphi(n)$ we get
\begin{equation}\label{3.4}
    |\varphi(n)\psi(n) - \varphi(m)\psi(m)|^2 \geq (\varphi(n)-\varphi(m))(\psi(n)^2 \varphi(n) - \psi(m)^2 \varphi(m)).
\end{equation}
Since $\varphi(n) \geq \varphi(n)-\varphi(m)$, the above inequality is true even when $\psi(m)=0$.
Using \eqref{3.4} and \eqref{3.1} we obtain
\begin{align*}
    \sum_{n=1}^\infty |u(n)-u(n-1)|^2v(n)&=
    \sum_{n=1}^\infty |\varphi(n)\psi(n)-\varphi(n-1)\psi(n-1)|^2 v(n)\\
    &\geq \sum_{n=1}^\infty \Big(\varphi(n) - \varphi(n-1)\Big)\Big(\psi(n)^2 \varphi(n) - \psi(n-1)^2 \varphi(n-1)\Big)v(n)\\
    & = \sum_{n=1}^\infty \Big(\frac{\Delta \varphi}{\varphi} v - \frac{(\varphi(n+1)-\varphi(n))(v(n+1)-v(n))}{\varphi}\Big)|u(n)|^2 \\
   & \geq \sum_{n=1}^\infty w(n) |u(n)|^2.
\end{align*}
This completes the proof.
\end{proof}
Now we are ready to prove Theorem \ref{thm2.1}.
\begin{proof}[Proof of theorem \ref{thm2.1}]
Let $v(n):= n^\alpha$ and $\varphi(n):= n^\beta $ on $\mathbb{N}$ and $\varphi(0):=0$ and $w_{\alpha, \beta}$ be as defined by \eqref{2.2}. It can be easily checked that the triplet $(v,\varphi, w)$ satisfies \eqref{3.1}. Now Theorem \ref{thm2.1} directly follows from the Lemma \ref{lem3.1}.
\end{proof}

In the next section we would be concerned about finding the parameters $\alpha$ and $\beta$ for which the weight $w_{\alpha, \beta}$ can be estimated from below by $\frac{(\alpha-1)^2}{4} n^{\alpha-2}$.

\section{Proof of Corollaries \ref{cor2.3} and \ref{cor2.5}}\label{sec4}
The goal in this section is to find parameters $\alpha$ and $\beta$ for which $w_{\alpha, \beta}(n) \geq \frac{(\alpha-1)^2}{4}n^{\alpha-2}$. With this in mind, we introduce the function $g_{\alpha, \beta}(x):= 1+(1+x)^\alpha - (1-x)^\beta - (1+x)^{\alpha+\beta}$. The goal now becomes to find parameters $\alpha$ and $\beta$ for which  
\begin{align*}
    g_{\alpha,\beta}(x) \geq \frac{(\alpha-1)^2}{4}x^2
\end{align*}
for $0 < x \leq 1/2$ and $w_{\alpha,\beta}(1)= 1+2^\alpha -2^{\alpha+\beta} \geq (\alpha-1)^2/4$.\\

Recall that, for $x \in (0,1)$, the Taylor series gives
\begin{equation}\label{4.1}
    (1\pm x)^r = \sum_{k=0}^\infty {r \choose k} (\pm1)^k x^k.
\end{equation}

Using \eqref{4.1}, we get the following expansion of $g_{\alpha,\beta}(x)$
\begin{equation}\label{4.2}
    g_{\alpha, \beta}(x) = \sum_{k=2}^\infty \Bigg[{\alpha \choose k} -(-1)^k {\beta \choose k} - {\alpha+\beta \choose k}\Bigg]x^k.
\end{equation}

Observe that the coefficient of $x^2$ is maximized when $ \beta = (1-\alpha)/2$. Taking $\beta = (1-\alpha)/2$, 
\begin{equation}\label{4.3}
    g(x) := g_{\alpha, \beta}(x) = \frac{(\alpha-1)^2}{4} x^2 + \sum_{k=3}^\infty \Bigg[{\alpha \choose k} -(-1)^k {(1-\alpha)/2 \choose k} - {(1+\alpha)/2 \choose k}\Bigg]x^k 
\end{equation}

In the next Lemma, we prove that the coefficients of $x^k$ in \eqref{4.3} are non-negative for $\alpha \in [1/3,1) \cup\{0\}$, which will be used as an ingredient in the proof of Corollary \ref{cor2.5}.
\begin{lem}\label{lem4.1} 
Let $b_k(\alpha)$ be defined as 
\begin{align*}
    b_k(\alpha) := {\alpha \choose k} -(-1)^k {(1-\alpha)/2 \choose k} - {(1+\alpha)/2 \choose k}.
\end{align*}
Then $b_k(\alpha) \geq 0$ for $\alpha \in [1/3,1) \cup\{0\}$ and $k \geq 3$.
\end{lem}

\begin{proof}
We will first cover the case $\alpha =0$. For $ k\geq 3$ we have
\begin{align*}
    b_k(0) = -(-1)^k {1/2 \choose k} - {1/2 \choose k} = -{1/2 \choose k}(1+ (-1)^k).
\end{align*}

Clearly, for odd $k$, $b_k(0) = 0$ and for even $k$ we have $b_k(0) = -2{1/2 \choose k}$, which is non-negative. This proves the non-negativity of $b_k(0)$ for $k \geq 3$.\\

Next we assume that $\alpha \in [1/3,1)$. Let $\alpha_1 := (1-\alpha)/2$ and $\alpha_2 := (1+\alpha)/2$. Then
\begin{align*}
    b_k(\alpha) &= {\alpha \choose k} - (-1)^k{\alpha_1 \choose k} -{\alpha_2 \choose k}\\
    &= (-1)^{k-1}\frac{\alpha(1-\alpha)...(k-1-\alpha)}{k!} + \frac{\alpha_1(1-\alpha_1)...(k-1 -\alpha_1)}{k!} + (-1)^k\frac{\alpha_2(1-\alpha_2)....(k-1-\alpha_2)}{k!} .
\end{align*}
We will treat the case of odd and even $k$ separately. First consider the case when $k$ is odd.
\begin{align*}
    b_k(\alpha) = {\alpha \choose k} + \frac{\alpha_1(1-\alpha_2)..(k-1-\alpha_2)}{k!}\Bigg[\prod_{i=1}^{k-1} \frac{(i-\alpha_1)}{(i-\alpha_2)} - \frac{\alpha_2}{\alpha_1}\Bigg] = {\alpha \choose k} + \frac{\alpha_1}{\alpha_2}{\alpha_2 \choose k}\Bigg[\prod_{i=1}^{k-1} \frac{(i-\alpha_1)}{(i-\alpha_2)} - \frac{\alpha_2}{\alpha_1}\Bigg]
\end{align*}
Note that for $i \geq 1$ we have $\frac{i-\alpha_1}{i-\alpha_2} = \frac{2i-1+\alpha}{2i-1-\alpha} \geq 1$. Therefore we have
\begin{align*}
    \prod_{i=1}^{k-1}\frac{(i-\alpha_1)}{(i-\alpha_2)} - \frac{\alpha_2}{\alpha_1} = \Big(\prod_{i=2}^{k-1}\frac{(i-\alpha_1)}{(i-\alpha_2)} - 1\Big) \frac{\alpha_2}{\alpha_1} \geq 0.
\end{align*}
The above inequality along with non-negativity of ${\alpha \choose k}, {\alpha_2 \choose k}$ for odd $k$ proves that, $b_k(\alpha) \geq 0$ for odd $k \geq 3$.\\

Next we consider the case when $k$ is even. 
\begin{equation}\label{4.4}
    \begin{split}
        b_k(\alpha) &= - \frac{\alpha(1-\alpha)...(k-1-\alpha)}{k!} + \frac{\alpha_1(1-\alpha_1)...(k-1 -\alpha_1)}{k!} - {\alpha_2 \choose k}\\
        &= \frac{\alpha_1(1-\alpha)...(k-1-\alpha)}{k!}\Big(\prod_{i=1}^{k-1} \frac{i-\alpha_1}{i-\alpha} - \frac{\alpha}{\alpha_1}   \Big) -  {\alpha_2 \choose k}.    
    \end{split}
\end{equation}

Consider the following polynomial in $\alpha$:
\begin{align*}
    P(\alpha) &:= \prod_{i=1}^{7} \frac{i-\alpha_1}{i-\alpha} - \frac{\alpha}{\alpha_1}\\
    &= \prod_{i=1}^{7} \frac{2i-1+\alpha}{2(i-\alpha)} - \frac{2\alpha}{1-\alpha} = \frac{1}{\prod_{i=1}^{7}2(i-\alpha)}Q(\alpha).
\end{align*}
where 
\begin{equation}\label{4.5}
    Q(\alpha):= \prod_{i=1}^7 (2i-1+\alpha) - 2^8 \alpha \prod_{i=2}^7(i-\alpha).
\end{equation}

Next we will show that $Q(\alpha)$ is non-negative for $\alpha \in [1/3,1)$. Note that showing $Q(\alpha) \geq 0$ is equivalent to showing 
\begin{equation}\label{4.6}
    \log(\prod_{i=1}^7 (2i-1+\alpha)) \geq \log(2^8 \alpha \prod_{i=2}^7(i-\alpha)).
\end{equation}

We introduce
\begin{align*}
    R(\alpha) &:= \log(\prod_{i=1}^7 (2i-1+\alpha)) - \log(2^8 \alpha \prod_{i=2}^7(i-\alpha))\\
    &= \sum_{i=1}^7 \log(2i-1+\alpha) - \log(2^8) - \log(\alpha) -\sum_{i=2}^7 \log(i-\alpha).
\end{align*}

It is straightforward to check that $R''(\alpha) \geq 0$ whenever $1/3 \leq \alpha \leq 1$. This, along with the fact that $R'(1/3)$ is non-negative, implies that $R'(\alpha) \geq 0$ in the specified domain. This means that the function $R(\alpha)$ is non-decreasing in the interval $(1/3, 1)$. Since $R(1/3) =0$, we can conclude that $R(\alpha) \geq 0$ in the interval $(1/3,1)$. Therefore we have $Q(\alpha) \geq 0$ which further implies that $P(\alpha)$ is non-negative in the interval $[1/3,1)$.\\

Also note that $\frac{i-\alpha_1}{i-\alpha} \geq 1$ for $1/3 \leq \alpha \leq 1$.  Using this fact along with the non-negativity of $P(\alpha)$ in \eqref{4.4} we get 
\begin{equation}\label{4.7}
    b_k(\alpha) \geq 0
\end{equation}\label{3.8}
for even $k\geq 8$ and $1/3 \leq \alpha <1$. \\

Now it remains to show that $b_4(\alpha)$ and $b_6(\alpha)$ are non-negative. Doing standard computations, we find that
\begin{equation}\label{4.8}
    b_4(\alpha) = \frac{1}{192}(5-\alpha)(1-\alpha)(7\alpha^2 - 6\alpha +3).
\end{equation}
and 
\begin{equation}\label{4.9}
    b_6(\alpha) = \frac{1}{23040}(1-\alpha)(9-\alpha)(31 \alpha^4 - 170 \alpha^3 + 536\alpha^2 - 310 \alpha + 105).
\end{equation}

It is very easy to see that $b_4(\alpha)$ is non-negative for $ 0 \leq \alpha < 1.$ Consider 
\begin{align*}
    T(\alpha) &:= 31 \alpha^4 - 170 \alpha^3 + 536\alpha^2 - 310 \alpha + 105.
\end{align*}

Let $\alpha^* := 7/20$. It can be easily verified that $T''(\alpha) \geq 0$ and both $T'(\alpha^*), T(\alpha^*)$ are non-negative. This implies the non-negativity of $T(\alpha)$ for $\alpha \in [\alpha^*,1)$. \\

Now assume $\alpha \in [0,\alpha^*]$. Using arithmetic-geometric mean inequality we get
\begin{align*}
    31\alpha^4 + 536\alpha^2 \geq 2\sqrt{16616}\alpha^3.
\end{align*}

Now showing $T(\alpha)$ is non-negative boils down to showing $\Tilde{T}(\alpha):=2\sqrt{16616}\alpha^3 -170 \alpha^3 -310 \alpha +105 \geq 0$. Observing that $\Tilde{T}'(\alpha) \leq 0$ for $\alpha \in (0,1)$ and $\Tilde{T}(\alpha^*) \geq 0$ proves the non-negativity of $\Tilde{T}(\alpha)$ in the interval $[0,\alpha^*]$. This proves the non-negativity of $T$ and hence the non-negativity of $b_6(\alpha)$ in the interval $ \alpha \in [0,1)$.   
\end{proof}

Next we will prove that $g(x) \geq \frac{(\alpha-1)^2}{4}x^2$ for $\alpha \in [0,1) \cup [5, \infty)$. We will treat the cases $\alpha \in [0,1)$ and when $\alpha \in [5,\infty)$ separately.
\begin{lem}\label{lem4.2}
Let $\alpha \in[0,1/3]$. Then 
\begin{equation}\label{4.10}
    g(x) \geq \frac{(\alpha-1)^2}{4} x^2
\end{equation}
for $0< x <1$.
\end{lem}

\begin{proof}
Let $E(x):= g(x) - \frac{(\alpha-1)^2}{4}x^2 =  1+(1+x)^\alpha - (1-x)^{(1-\alpha)/2} - (1+x)^{(1+\alpha)/2} -\frac{(\alpha-1)^2}{4}x^2$.\\

The first four derivatives of $E$ are given by
\begin{align*}
    E'(x) &= \alpha(1+x)^{\alpha-1} + \frac{1-\alpha}{2}(1-x)^{\frac{-1-\alpha}{2}} - \frac{(1+\alpha)}{2}(1+x)^{\frac{\alpha-1}{2}} - \frac{(\alpha-1)^2}{2}x.\\
    E''(x) &= \alpha(\alpha-1)(1+x)^{\alpha-2} + \frac{(1+\alpha)(1-\alpha)}{4}(1-x)^{\frac{-3-\alpha}{2}} + \frac{(1+\alpha)(1-\alpha)}{4}(1+x)^{\frac{\alpha-3}{2}} - \frac{(\alpha-1)^2}{2}.\\
    E'''(x) &= \alpha(\alpha-1)(\alpha-2)(1+x)^{\alpha-3} + \frac{(1+\alpha)(1-\alpha)(3+\alpha)}{8}(1-x)^{\frac{-5-\alpha}{2}}\\
    &+ \frac{(1+\alpha)(1-\alpha)(\alpha-3)}{8}(1+x)^{\frac{\alpha-5}{2}}.\\
    E''''(x) &= \alpha(\alpha-1)(\alpha-2)(\alpha-3)(1+x)^{\alpha-4} + \frac{(1+\alpha)(1-\alpha)(3+\alpha)(5+\alpha)}{16}(1-x)^{\frac{-7-\alpha}{2}}\\
    &+ \frac{(1+\alpha)(1-\alpha)(\alpha-3)(\alpha-5)}{16}(1+x)^{\frac{\alpha-7}{2}}. 
\end{align*}

Note that $E(0)=E'(0)=E''(0)=0$ and $E'''(0) = \frac{3}{4}\alpha(1-\alpha)(3-\alpha)$ which is non-negative. Further assuming that $E''''(x)$ is non-negative completes the proof. In what follows we will prove that $E''''(x)$ is non-negative. \\

Using arithmetic-geometric mean inequality we get
\begin{align*}
    2\Bigg( \frac{(1+\alpha)^2(9-\alpha^2)(25-\alpha^2)}{16^2}(1+x)^{\frac{\alpha-7}{2}}(1-x)^{\frac{-7-\alpha}{2}}\Bigg)^{\frac{1}{2}} &\leq \frac{(1+\alpha)(3+\alpha)(5+\alpha)}{16}(1-x)^{\frac{-7-\alpha}{2}}+\\
    &\frac{(1+\alpha)(3-\alpha)(5-\alpha)}{16}(1+x)^{\frac{\alpha-7}{2}}.
\end{align*}

Therefore proving $E''''(x)\geq 0$ reduces to showing 
\begin{equation}\label{4.11}
    2\Bigg(\frac{(1+\alpha)^2(9-\alpha^2)(25-\alpha^2)}{16^2}(1+x)^{\frac{\alpha-7}{2}}(1-x)^{\frac{-7-\alpha}{2}}\Bigg)^{\frac{1}{2}} \geq \alpha(2-\alpha)(3-\alpha)(1+x)^{\alpha-4}.
\end{equation}

which is equivalent to proving
\begin{align*}
    \log2 + 1/2 \log\Bigg(\frac{(1+\alpha)^2(9-\alpha^2)(25-\alpha^2)}{16^2}(1+x)^{\frac{\alpha-7}{2}}(1-x)^{\frac{-7-\alpha}{2}}\Bigg)\geq \log\Big(\alpha(2-\alpha)(3-\alpha)(1+x)^{\alpha-4}\Big).
\end{align*}

Consider the function 
\begin{align*}
    f(x) &:= \log2 + 1/2 \log\Bigg(\frac{(1+\alpha)^2(9-\alpha^2)(25-\alpha^2)}{16^2}(1+x)^{\frac{\alpha-7}{2}}(1-x)^{\frac{-7-\alpha}{2}}\Bigg)\\ &-\log\Big(\alpha(2-\alpha)(3-\alpha)(1+x)^{\alpha-4}\Big) \\
    &= \log2 + 1/2\log\Bigg(\frac{(1+\alpha)^2(9-\alpha^2)(25-\alpha^2)}{16^2}\Bigg) - \log\Big(\alpha(2-\alpha)(3-\alpha)\Big)\\
    &+ \frac{3}{4}(3-\alpha)\log(1+x) - \frac{7+\alpha}{4}\log(1-x).
\end{align*}

It can be easily checked that $f'(x)\geq 0$. Now we will show that $f(0)$ is non-negative for $\alpha \in (0,1/3]$. Consider
\begin{align*}
    2f(0)&= \log4  + \log\Bigg(\frac{(1+\alpha)^2(9-\alpha^2)(25-\alpha^2)}{16^2}\Bigg) - 2\log\Big(\alpha(2-\alpha)(3-\alpha)\Big)\\
    &= \log4 + \log\Bigg(\frac{(1+\alpha)^2(9-\alpha^2)(25-\alpha^2)}{16^2\alpha^2(2-\alpha)^2(3-\alpha)^2}\Bigg)\\
    &=\log4 + \log\Bigg(\frac{(1+\alpha)^2(3+\alpha)(25-\alpha^2)}{16^2\alpha^2(2-\alpha)^2(3-\alpha)}\Bigg).
\end{align*}

So $f(0)$ is non-negative iff 
\begin{equation}\label{4.12}
    \frac{(1+\alpha)^2(3+\alpha)(25-\alpha^2)}{16^2\alpha^2(2-\alpha)^2(3-\alpha)} \geq 1/4.
\end{equation}

Consider the function
\begin{align*}
    Q(\alpha) := (1+\alpha)^2(25-\alpha^2) - 64\alpha^2(2-\alpha)^2.
\end{align*}

It is straightforward to check that $Q''(\alpha)$ is negative in the interval $(0,1/3)$ and $Q'(0), Q(0)$ and  $Q(1/3)$ are non-negative. From this information one can easily conclude that $Q(\alpha) \geq 0$ in the interval $(0,1/3]$. Now consider
\begin{align*}
    \frac{(1+\alpha)^2(3+\alpha)(25-\alpha^2)}{16^2\alpha^2(2-\alpha)^2(3-\alpha)} \geq \frac{(1+\alpha)^2(25-\alpha^2)}{16^2\alpha^2(2-\alpha)^2} \geq 1/4.
\end{align*}

The last steps follows from the non-negativity of $Q(\alpha)$. This proves that $f(0)$ is non-negative whenever $\alpha \in (0,1/3]$. This fact, along with the non-negativity of $f'(x)$, implies $f(x) \geq 0$, which further implies $E''''(x) \geq 0$.
\end{proof}

\begin{rem}\label{rem3.3}
Using numerics, one can easily conclude that \eqref{4.12} is true for $\alpha \in (0,1)$. Therefore Lemma \ref{lem4.2} is true for $\alpha \in (0,1)$, i.e, $g(x) \geq \frac{(\alpha-1)^2}{4}x^2$ whenever $\alpha \in (0,1)$. But proving \eqref{4.12} in the interval $(0,1)$ mathematically becomes a bit tricky.
\end{rem}

\begin{rem}
Lemma \ref{lem4.1} along with Lemma \ref{lem4.2} proves that $g(x) \geq \frac{(\alpha-1)^2}{4}x^2$ for $x \in [0,1)$ and $\alpha \in  [0,1)$.
\end{rem}
Next we will prove that $g(x) \geq \frac{(\alpha-1)^2}{4}x^2$ when $\alpha \geq 5$.
\begin{lem}\label{lem4.5}
Let $\alpha \geq 5$. Then 
\begin{equation}\label{4.13}
    g(x) \geq \frac{(\alpha-1)^2}{4} x^2
\end{equation}
for $0 < x \leq 1/2$.
\end{lem}

\begin{proof}
Consider
\begin{align*}
    E(\alpha,x):=  1+(1+x)^{2\alpha+1}- (1-x)^{-\alpha}-(1+x)^{\alpha+1}-\alpha^2 x^2.
\end{align*}

Note that, under the transformation $\alpha \mapsto 2\alpha+1$, showing \eqref{4.13} reduces to proving $E(\alpha,x) \geq 0$ for $\alpha \geq 2$. The first three derivatives of $E$ w.r.t $\alpha$ are given by
\begin{align*}
    \partial_\alpha E(\alpha, x) &= 2(1+x)^{2\alpha+1}\log(1+x) +(1-x)^{-\alpha}\log(1-x) - (1+x)^{\alpha+1}\log(1+x) -2\alpha x^2.\\
    \partial^2_{\alpha^2}E(\alpha, x) &= 4(1+x)^{2\alpha+1}\log^2(1+x) - (1-x)^{-\alpha}\log^2(1-x) - (1+x)^{\alpha+1}\log^2(1+x) - 2x^2.\\
    \partial^3_{\alpha^3}E(\alpha,x) &= 8(1+x)^{2\alpha+1}\log^3(1+x) + (1-x)^{-\alpha}\log^3(1-x) - (1+x)^{\alpha+1}\log^3(1+x).
\end{align*}

The strategy of the proof is to show that $\partial^3_{\alpha^3}E(\alpha,x), \partial^2_{\alpha^2}E(2, x), \partial_\alpha E(2, x)$ and $E(2,x)$ are all non-negative, thereby completing the proof. \\

Consider
\begin{align*}
    \partial^3_{\alpha^3}E(\alpha,x) &= 8(1+x)^{2\alpha+1}\log^3(1+x) + (1-x)^{-\alpha}\log^3(1-x) - (1+x)^{\alpha+1}\log^3(1+x) \\
    &= [8(1+x)^{2\alpha+1}-(1+x)^{\alpha+1}]\log^3(1+x) + (1-x)^{-\alpha}\log^3(1-x)\\
    &= -(1-x)^{-\alpha}\log^3(1-x)\Big[(1+x)(1-x)^\alpha(8(1+x)^{2\alpha}-(1+x)^\alpha)\frac{\log^3(1+x)}{-\log^3(1-x)} - 1\Big]\\
    &= -(1-x)^{-\alpha}\log^3(1-x)\Big[(1+x)[8((1+x)^2(1-x))^\alpha-(1-x^2)^\alpha]\frac{\log^3(1+x)}{-\log^3(1-x)} - 1\Big]\\
    & \geq -(1-x)^{-\alpha}\log^3(1-x)\Big[(1+x)[8((1+x)^2(1-x))^2-(1-x^2)^2]\frac{\log^3(1+x)}{-\log^3(1-x)} - 1\Big]\\
    & = -(1-x)^{-\alpha}\log^3(1-x)\Big[(1+x)^3(1-x)^2[8(1+x)^2-1]\frac{\log^3(1+x)}{-\log^3(1-x)} - 1\Big]\\
    & \geq -(1-x)^{-\alpha}\log^3(1-x)\Big[7/4(1+x)^3\frac{\log^3(1+x)}{-\log^3(1-x)} - 1\Big].
\end{align*}

Therefore, for $\alpha \geq 2$, we have
\begin{equation}\label{4.14}
    \partial^3_{\alpha^3}E(\alpha,x)  \geq -(1-x)^{-\alpha}\log^3(1-x)\Big[7/4(1+x)^3\frac{\log^3(1+x)}{-\log^3(1-x)} - 1\Big].
\end{equation}

Next we will prove the following inequalities for $0<x \leq 1/2$.
\begin{align*}
    & 7/4(1+x)^3 \log^3(1+x) + \log^3(1-x) \geq 0.\\
    &\partial^2_{\alpha^2}E(2, x) = 4(1+x)^5\log^2(1+x) - (1-x)^{-2} \log^2(1-x) - (1+x)^3\log^2(1+x) - 2x^2 \geq 0.\\
    &\partial_\alpha E(2, x) = 2(1+x)^5\log(1+x) + (1+x)^{-2}\log(1-x) - (1+x)^3\log(1+x) - 4x^2 \geq 0. \\
    &E(2,x) = 1+(1+x)^5 - (1-x)^{-2}-(1+x)^3 - 4x^2 \geq 0.
\end{align*}

Assuming the above inequalities are true, the result follows.\\

Standard computations yield
\begin{align*}
    E_1(x) &:= \partial_\alpha E(2, x) =2(1+x)^5\log(1+x) + (1-x)^{-2}\log(1-x) - (1+x)^3\log(1+x) - 4x^2.\\
    E_1^{(5)}(x) &= 240\log(1+x) + \frac{6}{(1+x)^2} + 548 - \frac{1044}{(1-x)^7} + 720 \frac{\log(1-x)}{(1-x)^7}\\
    & \leq 240\log(3/2) + 6 + 548 - 1044 \leq 0
\end{align*}

It can be easily checked that $E_1^{(i)}(0) \geq 0$ for $i\leq 4$ and $E_1(1/2) \geq 0$. This proves that $E_1^{(i)}(x)$ for $1 \leq i \leq 4$ is either non-negative or it has one zero say $y$, such that $E_1^{(i)}(x) \geq 0$ for $ x \leq y$ and $E_1^{(i)}(x) \leq 0$ for $ x \geq y$. Let us assume that $E_1^{(1)}(x)$ is non-negative, this implies that $E_1(x)$ is a non-decreasing function of $x$. This combined with the fact that $E_1(0) = 0$ proves that $E_1(x) \geq 0$. Another possibility is that $E_1^{(1)}(x)$ has one zero $y$. Then $E_1(x)$ is a non-deceasing function in $[0,y]$ and it is non-increasing in $[y, 1/2]$. This combined with non-negativity of $E_1(1/2)$ proves that $E_1(x)\geq 0$ in the interval $(0,1/2]$.\\\\

Now consider the second derivative 
\begin{align*}
    E_2(x) &:= \partial^2_{\alpha^2}E(2, x) = 4(1+x)^5\log^2(1+x) - (1-x)^{-2} \log^2(1-x) - (1+x)^3\log^2(1+x) - 2x^2.\\
    E_2^{(6)}(x) &= \frac{3(40x^2 + 80x + 39)}{(1+x)^3}\log(1+x) + \frac{274}{1+x} - \frac{1}{(1+x)^3} -\frac{1276}{(1-x)^8}\\ &+\frac{9(223-70\log(1-x))}{(1-x)^8}\log(1-x) \\
    & \leq 267\log(3/2) + 274 - 8/27 - 1276\leq 0.
\end{align*}

Simple calculations yield $E_2^{(i)}(0) \geq 0$ for $i \leq 5$ and $E_2(1/2) \geq 0$. This proves that $E_2(x) \geq 0$ for $x \in (0,1/2]$, via the same logic used in proving that $E_1(x)$ is non-negative. \\

Next, we consider the third derivative
\begin{align*}
    E_3(x) &:= 7/4(1+x)^3 \log^3(1+x) + \log^3(1-x).\\
    E_3^{(5)}(x) & = - \frac{210}{(1-x)^5} + \frac{105}{2(1+x)^2} + 300 \frac{\log(1-x)}{(1-x)^5} - 72 \frac{\log^2(1-x)}{(1-x)^5} - \frac{105}{2} \frac{\log(1+x)}{(1+x)^2} - \frac{63}{2} \frac{\log^2(1+x)}{(1+x)^2} \\
    & \leq -210 + 105/2 \leq 0.
\end{align*}

Furthermore, $E_3^{(i)}(0) \geq 0$ for $i \leq 4$ and $E_3(1/2)\geq 0$. This proves the non-negativity of $E_3(x)$. \\

Finally, we consider $E(2,x)$
\begin{align*}
    E_0(x) &:= E(2,x) = 1+(1+x)^5 - (1-x)^{-2}-(1+x)^3 - 4x^2.\\
    E_0^{(5)}(x) &= 120 -\frac{720}{(1-x)^7} \leq 0.
\end{align*}

It can be verified that $E_0^{(i)}(0) \geq 0$ for $i \leq 4$ and $E_0(1/2) \geq 0$. This implies that $E_0(x) \geq 0$ in the interval $(0,1/2]$. 
\end{proof}

\begin{rem}
Using Lemmas \ref{lem4.1}, \ref{lem4.2} and \ref{lem4.5} we can conclude that $g(x) \geq \frac{(\alpha-1)^2}{4}x^2$ for $0 < x \leq 1/2$ and $\alpha \in [0,1) \cup [5,\infty)$. This proves that, with the choice $\beta = (1-\alpha)/2$, we have $w_{\alpha,\beta}(n) \geq \frac{(\alpha-1)^2}{4}n^{\alpha-2}$ for $n \geq 2$ and $\alpha \in [0,1) \cup [5,\infty)$. Now it remains to show that $w_{\alpha, (1-\alpha)/2}(1) \geq (\alpha-1)^2/4$. This will be proved in the next Lemma.
\end{rem}

\begin{lem}\label{lem4.7}
Let $w_{\alpha, \beta}$ be the weight function as defined by \eqref{2.2}. Then for $\beta= (1-\alpha)/2$ and  $\alpha \in [0,1) \cup [5, \infty)$ we have
\begin{equation}\label{4.15}
    w_{\alpha,\beta}(1) = 1+2^\alpha - 2^{(\alpha+\beta)/2}=1+2^\alpha -2^{(1+\alpha)/2} \geq \frac{(\alpha-1)^2}{4}.
\end{equation}
\end{lem}

\begin{proof}
We will consider the case, when $\alpha \in [0,1)$ and $\alpha \geq 5$ separately. First assume $\alpha \in [0,1)$. Using mean value theorem for the function $2^x$, we get, for $\xi \in [\alpha, (1+\alpha)/2]$, 
\begin{align*}
    2^{(1+\alpha)/2} - 2^\alpha = \frac{(1-\alpha)}{2} 2^\xi \log2 \leq \frac{(1-\alpha)}{2} 2^{(1+\alpha)/2}\log2.
\end{align*}

This implies that 
\begin{align*}
    w_{\alpha, (1-\alpha)/2}(1)  - \frac{(\alpha-1)^2}{4} \geq 1 - 2^{(1+\alpha)/2}\log2\frac{(1-\alpha)}{2} - \frac{(\alpha-1)^2}{4} =: g(\alpha).
\end{align*}

Derivatives of $g$ are given by
\begin{align*}
    g'(\alpha) &= 1/2[2^{(\alpha+1)/2}\log2 - 2^{(\alpha+1)/2}\frac{(1-\alpha)}{2}\log^22 - \alpha+1].\\
    g''(\alpha) &= 1/4[2^{(\alpha+3)/2}\log^22 - 2^{(\alpha+1)/2}\frac{1-\alpha}{2}\log^32 -2].\\
    g'''(\alpha) &= \frac{2^{(1+\alpha)/2}}{8}\log^32[3-(1-\alpha)/2 \log2] \geq 0.  
\end{align*}

Note that $g''(1) = \log^2(2)-1/2 \leq 0$, $g'(1) = \log2 \geq 0$ and $g(0) = (3-2\sqrt{2}\log(2))/4 \geq 0$. From this we can conclude that $w_{\alpha,(1-\alpha)/2}(1) \geq \frac{(\alpha-1)^2}{4}$ for $\alpha \in [0,1)$.\\

Now let $\alpha \geq 5$ case. Let $h(\alpha):= 1+ 2^{\alpha} - 2^{(1+\alpha)/2} - \frac{(\alpha-1)^2}{4}$. Derivatives of $h$ are given by
\begin{align*}
    h'(\alpha) &= 2^\alpha \log2 - \frac{2^{(1+\alpha)/2}}{2}\log2 - \frac{(\alpha-1)}{2}.\\
    h''(\alpha) &= 2^\alpha \log^22 - \frac{2^{(1+\alpha)/2}}{4}\log^22 - 1/2.\\
    h'''(\alpha) &= 2^\alpha \log^32 - \frac{2^{(1+\alpha)/2}}{8}\log^32 = \log^32(2^\alpha -2^{(\alpha-5)/2}) \geq 0.
\end{align*}

Noting that $h''(5) = 30\log^32 -1/2 \geq 0$, $h'(5) = (28 \log2-2)\geq 0$ and $h(5) = 21 \geq 0$. This proves that $h(\alpha) \geq 0$ for $\alpha \geq 5$.   
\end{proof}

Now we have all the pieces required to prove the Corollaries \ref{cor2.3} and \ref{cor2.5}. Let us put them together.
\begin{proof}[Proof of Corollary \ref{cor2.3}] Using Lemma \ref{lem4.1}, Lemma \ref{lem4.2} and Lemma \ref{lem4.5} we can conclude that
\begin{equation}\label{4.16}
    g(x) = 1+(1+x)^\alpha -(1-x)^{(1-\alpha)/2}-(1+x)^{(1+\alpha)/2} \geq \frac{(\alpha-1)^2}{4}x^2
\end{equation}
for $0 < x \leq 1/2$ and $\alpha \in [0,1) \cup [5, \infty)$. Now taking $x=1/n$, we get, for $n \geq 2$,
\begin{equation}\label{4.17}
    1 + \Big(1+\frac{1}{n}\Big)^\alpha - \Big(1-\frac{1}{n}\Big)^{(1-\alpha)/2} - \Big(1+\frac{1}{n}\Big)^{(1+\alpha)/2} \geq \frac{(\alpha-1)^2}{4}\frac{1}{n^2}.
\end{equation}

Using \eqref{4.17} along with Lemma \ref{lem4.7}, we conclude that, for $\beta =(1-\alpha)/2$,
\begin{equation}\label{4.18}
    w_{\alpha, \beta}(n) \geq \frac{(\alpha-1)^2}{4} n^{\alpha-2}
\end{equation}
for all $n \geq 1$.\\

Inequality \eqref{4.18} along with Theorem \ref{thm2.1}(with $\beta = (1-\alpha)/2$) proves Corollary \ref{cor2.3}. Next we prove the sharpness of the constant in Corollary \ref{cor2.3}.\\

Let $C$ be a constant such that
\begin{equation}\label{4.19}
    \sum_{n=1}^\infty |u(n)-u(n-1)|^2 n^\alpha \geq C \sum_{n=1}^\infty |u(n)|^2 n^{\alpha-2}
\end{equation}
for all $u \in C_c(\mathbb{N}_0)$ and $u(0)=0$.\\

Let $N \in \mathbb{N}$, $\beta \in \mathbb{R}$  and $\alpha\geq 0$ such that $2\beta +\alpha-2 <-1$, in particular, $\beta < 1/2$. Consider the following family of finitely supported functions on $\mathbb{N}_0$. 
\begin{align*}
u_{\beta,N}(n):=
\begin{cases}
n^\beta \hspace{73pt} &\text{for} \hspace{5pt}  1 \leq  n \leq N \\
-N^{\beta-1} n + 2N^{\beta} \hspace{13pt}  &\text{for} \hspace{5pt}  N \leq n \leq 2N \\
0 \hspace{83pt} &\text{for} \hspace{5pt}  n \geq 2N \hspace{5pt} \text{and} \hspace{5pt} n=0
\end{cases} 
\end{align*}

Clearly we have
\begin{equation}\label{4.20}
    \sum_{n=1}^\infty|u_{\beta,N}(n)|^2n^{\alpha-2} \geq \sum_{n=1}^{N}n^{2\beta+\alpha-2}.
\end{equation}
and
\begin{equation}\label{3.21}
    \begin{split}
        \sum_{n=1}^\infty|u_{\beta, N}(n)-u_{\beta,N}(n-1)|^2n^{\alpha} &= \sum_{n=2}^{N}(n^\beta - (n-1)^\beta)^2n^{\alpha} + \sum_{n=N+1}^{2N}N^{2\beta-2}n^{\alpha} + 1.
    \end{split}
\end{equation}
Using the fact that $\beta <1/2$, we get the following  basic estimates:\\
$(n^\beta - (n-1)^\beta)^2 \leq \beta^2(n-1)^{2\beta-2}$.\\
$\sum_{n=N+1}^{2N} n^\alpha \leq \int_{N+1}^{2N+1} x^\alpha dx = \frac{(2N+1)^{\alpha+1} - (N+1)^{\alpha+1}}{\alpha+1}$.\\\\

Using the above, in \eqref{3.21}, we get
\begin{equation}\label{4.22}
    \begin{split}
        \sum_{n=1}^\infty|u_{\beta, N}(n)-u_{\beta, N}(n-1)|^2 n^\alpha &\leq \beta^2 \sum_{n=2}^{N}(n-1)^{2\beta-2}n^\alpha + 
        \frac{N^{2\beta+\alpha-1}}{\alpha+1}\Bigg[\Big(2+\frac{1}{N}\Big)^{\alpha+1} - \Big(1+\frac{1}{N}\Big)^{\alpha+1}\Bigg]
        + 1.
    \end{split}
\end{equation}

Using estimates \eqref{4.20} and \eqref{4.22} in \eqref{4.19}, and taking limit $N \rightarrow \infty$, we get
\begin{equation}\label{4.23}
    C\sum_{n=1}^{\infty}n^{2\beta+\alpha-2} \leq \beta^2\sum_{n=2}^{\infty}(n-1)^{2\beta-2}n^\alpha + 1.
\end{equation}

Using Taylor's theorem for the function $x^\alpha$, we get, for $n \geq 2$,
\begin{equation}\label{4.24}
    n^\alpha = (1+n-1)^\alpha \leq (n-1)^\alpha + {\alpha \choose 1}(n-1)^{\alpha-1} + .... + {\alpha \choose \ceil*{\alpha}}(n-1)^{\alpha - \ceil*{\alpha}}.  
\end{equation}
where $\ceil*{\alpha}$ denotes the smallest integer greater than or equal to $\alpha$.\\ 
Using \eqref{4.24} in \eqref{4.23}, we obtain
\begin{equation}\label{4.25}
    C\sum_{n=1}^{\infty}n^{2\beta+\alpha-2} 
    \leq \beta^2 \sum_{i=0}^{\ceil*{\alpha}}{\alpha \choose i}\sum_{n=1}^\infty n^{2\beta + \alpha - i -2} + 1. 
\end{equation}

Finally, taking limit $\beta \rightarrow \frac{1-\alpha}{2}$, and observing that $\text{lim sup}_{\beta \rightarrow(1-\alpha)/2} \sum_{n=1}^ \infty n^{2\beta + \alpha -i-2}$ is finite for $i \geq 1$ and is infinite for $i=0$, we obtain
\begin{equation}\label{4.26}
    C \leq \frac{(\alpha-1)^2}{4}. 
\end{equation}
\end{proof}

\begin{proof}[Proof of Corollary \ref{cor2.5}]
Let $g(x)$ be as defined by \eqref{4.3}, that is, $g(x) := 1 + (1+x)^\alpha - (1-x)^\beta - (1+x)^{\alpha+\beta}$ for $\beta = (1-\alpha/2)$. Using Taylor's expansion of $g(x)$ we get identity $\eqref{4.3}$ for $ x \in (0,1)$ 
\begin{equation}\label{4.27}
    g(x) = \frac{(\alpha-1)^2}{4}x^2 + \sum_{k=3}^\infty b_k(\alpha) x^k 
\end{equation}
where 
\begin{align*}
    b_k(\alpha) := {\alpha \choose k} -(-1)^k {(1-\alpha)/2 \choose k} - {(1+\alpha)/2 \choose k}.
\end{align*}
Taking $x=1/n$ and multiplying both sides of \eqref{4.27} by a factor of $n^\alpha$, we obtain
\begin{equation}\label{4.28}
    w_{\alpha, \beta}(n) = \frac{(\alpha-1)^2}{4}\frac{n^\alpha}{n^2} + \sum_{k=3}^\infty b_k(\alpha) \frac{n^\alpha}{n^k}
\end{equation}
for $\beta = (1-\alpha)/2$ and $n \geq 2$. Using \eqref{4.28} along with Lemma \ref{lem4.7} in Theorem \ref{thm2.1}(with $\beta = (1-\alpha)/2$) proves inequality \eqref{1.7} for $\alpha \in [0,1) \cup [5, \infty)$. Finally using Lemma \ref{lem4.1} to note the non-negativity of $b_k(\alpha)$ for $\alpha \in [1/3,1) \cup \{0\}$ we complete the proof of Corollary \ref{cor2.5}.
\end{proof}

\section{Limitations of the Method}\label{sec5}
In this section our first goal is to point out that the method described in this paper doesn't work for proving Corollary \ref{cor2.3} when $\alpha < 0$ or $\alpha \in (1,4)$. This will be proved in Lemma's \ref{lem5.1} and \ref{lem5.2}. Our second goal is to show that Corollary \ref{cor2.3} cannot be improved in the sense of Corollary \ref{cor2.5} when $\alpha$ doesn't lie in the interval $[1/3,1)$. This will be achieved partially via Lemma \ref{lem5.5}.

\begin{lem}\label{lem5.1}
Let $\alpha < 0$. Then $\exists$ $\epsilon >0$(depending on $\alpha$) such that $g(x) < \frac{(\alpha-1)^2}{4}x^2$ for all $x \in (0,\epsilon)$.
\end{lem}

\begin{proof}
Let $E(x):= g(x) - \frac{(\alpha-1)^2}{4}x^2$. Computations done in Lemma \ref{lem4.2} give $E(0) = E'(0) = E''(0) = 0$ and $E'''(0) = \frac{3}{4}\alpha(1-\alpha)(3-\alpha)$. Clearly $E'''(0) < 0$ for negative $\alpha$. The result now follows from the continuity of derivatives of $E(x)$.
\end{proof}

\begin{lem}\label{lem5.2}
Let $\alpha \in (1,4)$ then $\exists$ $\epsilon >0$(depending on $\alpha$) such that $g(x) < \frac{(\alpha-1)^2}{4}$ for all $x \in (1/2-\epsilon, 1/2)$.
\end{lem}

\begin{proof}
Let $E(x):= g(x) - \frac{(\alpha-1)^2}{4}x^2$. We show that $E(1/2)$ is negative whenever $\alpha \in (1,4)$. The result then follows from the continuity of the function $E(x)$.\\
Standard computations yield 
\begin{align*}
    f(\alpha) &:= E(1/2) = 1+(3/2)^\alpha - (1/2)^{(1-\alpha)/2} - (3/2)^{(1+\alpha)/2} - \frac{(\alpha-1)^2}{16}.\\
    f'(\alpha) &= (3/2)^\alpha \log(3/2) + \frac{1}{2}(1/2)^{(1-\alpha)/2}\log(1/2) - \frac{1}{2}(3/2)^{(1+\alpha)/2}\log(3/2) - \frac{\alpha-1}{8}.\\
    f''(\alpha) &= (3/2)^\alpha \log^2(3/2) - \frac{1}{4}(1/2)^{(1-\alpha)/2}\log^2(1/2)- \frac{1}{4}(3/2)^{(1+\alpha)/2}\log^2(3/2) -1/8.\\
    f'''(\alpha) &= (3/2)^\alpha \log^3(3/2) + \frac{1}{8} (1/2)^{(1-\alpha)/2}\log^3(1/2) - \frac{1}{8}(3/2)^{(1+\alpha)/2}\log^3(3/2)\\
    & = \frac{2^{(\alpha-1)/2}}{8}\log^3(2)\Bigg[\Big(8 \sqrt{2}\Big(\frac{3}{2\sqrt{2}}\Big)^\alpha - \sqrt{3}\Big(\frac{\sqrt{3}}{2}\Big)^\alpha\Big)\frac{\log^3(3/2)}{\log^3(2)}-1\Bigg]\\
    & \geq 2^{(\alpha-1)/2}\log^3(2)\Big(\frac{21}{2} \frac{\log^3(3/2)}{\log^3(2)}-1\Big)\geq 0.
\end{align*}
It can be easily seen that $f(4), f'(1), f''(1)$ are negative. Since $f''(1)$ is negative and $f'''(\alpha) \geq 0$, there are two possibilities. Firstly $f''(\alpha) \leq 0$. Then negativity of $f'(1)$ implies that $f'(\alpha)<0$, which further imply that $f(\alpha)$ is a strictly decreasing function. This along with $f(1) =0$ proves that $f(\alpha)<0$ for $\alpha \in (1,4)$.  Second possibility is that there exists $\beta \in (1,4)$ such that $f''(\alpha) < 0$ for $\alpha \in [1, \beta)$ and  $f''(\alpha) \geq 0$ for $\alpha \in [\beta, 4)$. Now we further have two possibilities, first $f'(\alpha) <0$, this along with $f(1)=0$ would imply that $f(\alpha)<0$ for $\alpha \in (1,4)$. Second possibility is that there exists $\gamma \in (1,4)$ such that $f'(\alpha)<0$ in $[1, \gamma)$ and $f'(\alpha) \geq 0$ in $[\gamma, 4)$. This along with $f(1) =0$ and $f(4) < 0$ implies that $f(\alpha) < 0$ for $\alpha \in (1,4)$.
\end{proof}

\begin{rem}\label{rem5.3}
Lemmas \ref{lem5.1} and \ref{lem5.2} say that the weights $w_{\alpha, \beta}$(with $\beta = (1-\alpha)/2$) obtained in Theorem $\ref{thm2.1}$ do not control the weight $\frac{(\alpha-1)^2}{4}n^{\alpha-2}$ whenever $\alpha < 0$ or $\alpha \in  (1,4)$. Therefore, one cannot obtain Corollary $\ref{cor2.3}$ from the Theorem \ref{thm2.1} when $\alpha < 0$ or $\alpha \in (1,4)$.
\end{rem}

\begin{rem}\label{rem5.4}
Using Theorem \ref{thm2.1}(with $\beta = (1-\alpha)/2$) and the Taylor expansion of $g$ \eqref{4.3}, and Lemma \eqref{lem4.7}, we conclude that \eqref{2.4} holds true for $\alpha \in [0,1) \cup [5,\infty)$. We conjecture that constants $b_k(\alpha)$ given by \eqref{1.8} are not non-negative for all $k$ when $\alpha$ does not lie in $[1/3,1)\cup \{0\}$, i.e, for all $\alpha \in (0,1/3) \cup [5, \infty)$, there exists $i \geq 1$ such that $b_i(\alpha) < 0$. Therefore, we don't have the improvement \eqref{2.4} of inequality \eqref{2.3} when $\alpha$ does not lie in $[1/3,1) \cup \{0\}$. In the next Lemma, we prove a result which supports the conjecture. 
\end{rem}

\begin{lem}\label{lem5.5}
Let  $b_i(\alpha)$ be as defined by \eqref{2.5}. Let $\alpha = 2k+1$ then we have
\begin{align*}
    b_i(2k+1) = {2k+1 \choose i} - (-1)^i {-k \choose i} - {k+1 \choose i} 
\end{align*}
If $k \geq 2$, then  
\begin{align*}
    b_i(2k+1) &\geq 0 \hspace{19pt} \text{for} \hspace{5pt} 2 \leq i \leq k+1\\
    b_i(2k+1) &< 0 \hspace{19pt} \text{for} \hspace{5pt} i > k+1
\end{align*}

\end{lem}

\begin{proof}
Clearly, ${2k+1 \choose i} = {k+1 \choose i} = 0$ for $i \geq 2k+2$. Therefore  $b_i(2k+1) < 0$ for $i \geq 2k+2$. \\

Consider $ k+1 < i \leq 2k+1 $. In this case, we have
\begin{align*}
    b_i(2k+1) &= {2k+1 \choose i} - (-1)^i {-k \choose i}\\
    &= \frac{1}{i!}\Big((2k+1)2k(2k-1)..(2k+1-(i-1)) - k(k+1)(k+2)...(k+i-1) \Big)\\
    &= \frac{1}{i!}k(k+1)...(2k+1)\Big((k-1)..(2k+1-(i-1)) - (2k+2)...(k+i-1) \Big) <  0.
\end{align*}

In the case when $2 \leq  i \leq k+1$, we have
\begin{align*}
    b_i(2k+1) &= {2k+1 \choose i} - (-1)^i {-k \choose i} - {k+1 \choose i} \\
    & = \frac{1}{i!} \Big((2k+1)2k...(2k+1-(i-1)) - k(k+1)
    ..(k+i-1) - (k+1)k...(k+1-(i-1))\Big)\\
    &\geq \frac{1}{i!} \Big((2k+1)2k...(2k+1-(i-1)) - 2k(k+1)
    ..(k+i-1)\Big)
\end{align*}

Observing that $(2k-1)..(2k+1-(i-1)) \geq (k+2)...(k+i-1)$ for $k \geq 3$ and $i \leq k$, we get $b_i(k) \geq 0$. Now consider $i =k+1$. 
\begin{align*}
    b_i(2k+1) &= \frac{1}{i!} \Big((2k+1)2k..(k+1) - k(k+1)..(2k) - (k+1)!\Big)\\
    & = \frac{1}{i!}\Big((k+1)(k+1)(k+2)..2k - (k+1)k....1\Big)\\
    &\geq 0.
\end{align*}

The only case that remains is when $k=2$ and $i= 2$. It is straightforward that $b_2(5) = 4 \geq 0$.
\end{proof}

\end{document}